\documentclass[smallextended,runningheads]{svjour3}
\journalname{Annals of the Institute of Statistical Mathematics}
\smartqed  
\usepackage[dvips]{graphicx}
\usepackage{amssymb}
\usepackage{amsmath}
\usepackage{natbib}
\usepackage{rotating}
\usepackage{lscape}
\usepackage{float}
\newcommand{\numgroup}{J}
\newcommand{\numsubgroup}{m}
\newcommand{\numjkn}{m}
\newcommand{\numjkc}{c}
\newcommand{\totalitems}{\tau}
\newcommand{\samplesize}{n}
\newcommand{\numrowsA}{d}
\newcommand{\numcolsA}{\nu}
\newcommand{\marginalprobs}{q}

\def\AAA{{{
A}_{\totalitems, {\bf b}, {\bf c}, {\bf r}, {\bf s}}}}
\def\Ac{{ 
A}}

\renewcommand{\tabcolsep}{4pt}

\usepackage{bm}
\bmdefine{\Bzero}{0}
\bmdefine{\Bone}{1}

\def\Bt{{\bf t}}
\def\Bx{{\bf x}}

\def\Bz{{\bf z}}

\def\Bi{{\bf i}}
\def\Bj{{\bf j}}

\def\Bp{{\bf p}}

\def\Bone{{\bf 1}}

\def\ab{{\bf a}}

\def\qb{{\bf q}}
\def\zb{{\bf z}}

\def\RR{{\mathbb R}}
\def\ZZ{{\mathbb Z}}
\def\NN{{\mathbb N}}

\begin{document}

\title{Markov basis and Gr\"obner basis of
Segre-Veronese configuration
for testing independence in group-wise selections}
\author{Satoshi Aoki 
\and Takayuki Hibi 
\and Hidefumi Ohsugi
\and Akimichi Takemura
}

\institute{
Satoshi Aoki \at
Department of Mathematics and Computer Science\\
Kagoshima University\\
1-21-35, Korimoto, Kagoshima, Kagoshima 890-0065, Japan\\
\email{aoki@sci.kagoshima-u.ac.jp}
\and
Takayuki Hibi \at
Graduate School of Information Science and Technology\\
Osaka University\\
1-1, Yamadaoka, Suita, Osaka 565-0871, Japan
\and
Hidefumi Ohsugi \at
Department of Mathematics\\
Rikkyo University\\
3-34-1, Nishi Ikebukuro, Toshima-ku, Tokyo 171-8501, Japan
\and
Akimichi Takemura \at
Graduate School of Information Science and Technology\\
University of Tokyo\\
7-3-1, Hongo, Bunkyo-ku, Tokyo 113-0033, Japan
}

\date{Received: date / Revised: date}

\titlerunning{Gr\"obner basis for testing independence in group-wise selections}
\maketitle

\begin{abstract}
We consider testing independence 
in group-wise selections with some
restrictions on combinations of choices.  We present models 
for frequency data of selections for which 
it is easy to perform conditional tests 
by Markov chain Monte Carlo (MCMC) methods.  When the 
restrictions on the combinations can be described in terms of
a Segre-Veronese configuration, an explicit form of a
Gr\"obner basis consisting of binomials of degree two 
is readily available for performing a Markov chain.
We illustrate our setting with the National Center Test
for university entrance examinations in Japan.  We also apply our
method to testing independence hypotheses 
involving  genotypes at more than one locus
or haplotypes of alleles on the same chromosome.
\end{abstract}

\section{Introduction}
Suppose that people are asked to select items  which are
classified into categories or groups and there are some
restrictions on combinations of choices.  For example, when a consumer
buys a car, he or she can choose various options, such as a color, a
grade of air conditioning, a brand of audio equipment, etc.  Due to space
restrictions for example, some combinations of options  may not be
available.  The problem we consider in this paper is testing
independence of people's preferences in group-wise selections in the
presence of restrictions.  We assume that observations are the
counts of people choosing various combinations in group-wise
selections, i.e., the data are given in a form of a multiway
contingency table with some structural zeros corresponding to the
restrictions.

If there are $m$ groups of items and a consumer freely chooses just
one item from each group, then the combination of choices is simply a
cell of an $m$-way contingency table.  Then the hypothesis of
independence reduces to the complete independence model of an $m$-way
contingency table.  The problem becomes harder if there are some
additional conditions in a group-wise selection.  A consumer may be
asked to choose up to two items from a group or there may be a
restriction on the total number of items.  Groups may be nested, so
that there are further restrictions on the number of items from
subgroups.  Some restrictions may concern several groups or subgroups.
Therefore the restrictions on combinations may be complicated.  

As a concrete example we consider restrictions on choosing subjects in
the {\em National Center Test} (NCT hereafter) for university entrance
examinations in Japan. Due to time
constraints of the schedule of the test, the pattern of restrictions
is rather complicated.  However we will show that restrictions of NCT
can be described in terms of a Segre-Veronese configuration.

Another important application of this paper is a generalization of the
Hardy-Weinberg model in population genetics. We are interested in
testing various hypotheses of independence involving genotypes at more
than one locus and haplotypes of combination of alleles on the same
chromosome.  Although this problem seems to be different from the
above introductory motivation on consumer choices, we can imagine
that each offspring is required to choose two alleles for each gene
(locus) from a pool of alleles for the gene.  He or she can choose the
same allele twice (homozygote) or different alleles (heterozygote). In
the Hardy-Weinberg model two choices are assumed to be independently
and identically distributed.  A natural generalization of the
Hardy-Weinberg model for a single locus is to consider independence of
genotypes of more than one locus.  In many epidemiological studies,
the primary interest is the correlation between a certain disease and
the genotype of a single gene (or the genotypes at more than one
locus, or the haplotypes involving alleles on the same chromosome).
Further complication might arise if certain homozygotes are fatal and
can not be observed, thus becoming a structural zero.

In this paper we consider conditional tests of independence hypotheses
in the above two important problems from the viewpoint of Markov bases and 
Gr\"obner bases.  
Evaluation of $P$-values by Markov chain Monte Carlo (MCMC) method using
Markov bases and Gr\"obner bases was initiated in Diaconis and
Sturmfels (1998).
See also Sturmfels (1995).
Since then, this approach attracted much
attention from statisticians as well as algebraists.
Contributions of the present authors are found, for example, in
Aoki and Takemura (2005, 2007), Ohsugi and Hibi (2005, 2006, 2007), 
and Takemura and Aoki (2004). 
Methods of algebraic statistics are currently actively applied to
problems in computational biology (Pachter and Sturmfels, 2005). 
In algebraic statistics, results in commutative algebra may find 
somewhat unexpected applications in statistics.  At the same time statistical
problems may present new problems to commutative algebra.
A recent example is a conjunctive Bayesian network proposed in 
Beerenwinkel {\it et al.} (2006), 
where a result of Hibi (1987) is successfully
used.  In this paper we present application of results on
Segre-Veronese configuration to testing independence in NCT and
Hardy-Weinberg models.  In fact, these statistical considerations have
prompted  further theoretical developments of Gr\"obner bases for 
Segre-Veronese type configurations
and we will present 
these theoretical results  in our subsequent paper (Aoki et al., 2007).

Even in two-way tables, if the positions of the structural zeros are
arbitrary, then Markov bases may contain moves of high degrees 
(Aoki and Takemura, 2005). 
See also Huber {\it et al.} (2006) and Rapallo (2006) for Markov bases
of the problems with the structural zeros.
However if the restrictions on the combinations
can be described in terms of a Segre-Veronese configuration, then an
explicit form of a Gr\"obner basis consisting of binomials of
degree two with a squarefree initial term 
is readily available for running a Markov chain for
performing conditional tests of various hypotheses of independence.
Therefore models which can be described by a Segre-Veronese
configuration are very useful for statistical analysis.

The organization of this paper is as follows.  In Section 2,
we introduce two examples of group-wise selection.
In Section 3, we give a formalization of conditional tests and
MCMC procedures and consider various hypotheses
of independence for NCT data and the allele frequency data.
In Section 4, we define Segre-Veronese configuration. We
give an explicit expression of a reduced Gr\"obner basis for the
configuration and describe a simple procedure for running MCMC using
the basis for conditional tests.
In Section 5 we present numerical results on NCT data and diplotype
frequencies data.
We end the paper by some discussions in Section 6.

\section{Examples of group-wise selections}
In this section, we introduce two examples of group-wise selection.
In Section 2.1, we take a close look at patterns of selections of
subjects in NCT.  In Section 2.2, we illustrate an important
problem of population genetics from the viewpoint of group-wise selection.

\subsection{The case of National Center Test in Japan}
One important example of group-wise selection is the entrance
examination for universities in Japan.  In Japan, as the common
first-stage screening process, most students applying for universities
take the National Center Test for university entrance
examinations administered by National Center for University Entrance Examinations
(NCUEE).  
Basic information in English on NCT in 2006 
is available from the booklet published by NCUEE ([12] in
the references).
After obtaining the score of NCT, students apply to departments of individual
universities and take second-stage examinations administered by the
universities.  Due to time constraints of the schedule of NCT, there are
rather complicated restrictions on possible combination of subjects.
Furthermore each department of each university can impose different 
additional requirement on the combinations of subjects of NCT to students
applying to the department.  

In NCT examinees can choose subjects in Mathematics, Social Studies
and Science.  These three major subjects are divided into
subcategories.  For example Mathematics is divided into Mathematics 1
and Mathematics 2 and these are then composed of individual subjects.
In the test carried out in 2006, examinees could select two
mathematics subjects, two social studies subjects and three science subjects
at most as shown below.  The details of the subjects can be found in
web pages and publications of NCUEE.  
In this paper, we omit Mathematics for simplicity, and only consider
selections in Social Studies and Science.
In parentheses we show our
abbreviations for the subjects in this paper.
\begin{itemize}
 \item Social Studies:
 \begin{itemize}
  \item[$\circ$] Geography and History:\ One subject from \{World History
		 A (WHA),\ World History B (WHB),\ Japanese History A
		 (JHA),\  Japanese History B (JHB),\ Geography A (GeoA),\ 
		 Geography B (GeoB)\}
  \item[$\circ$] Civics:\ One subject from \{Contemporary Society
		 (ContSoc),\ Ethics,\
		 Politics and Economics (P\&E)\}
 \end{itemize}
 \item Science:
 \begin{itemize}
  \item[$\circ$] Science $1$:\ One subject from \{Comprehensive Science
		 B (CSciB),\
		 Biology I (BioI),\ Integrated Science (IntegS),\
		 Biology IA (BioIA)\}
  \item[$\circ$] Science $2$:\ One subject from \{Comprehensive Science
		 A (CSciA),\
		 Chemistry I (ChemI),\ Chemistry IA (ChemIA)\}
  \item[$\circ$] Science $3$:\ One subject from \{Physics I (PhysI), 
		 Earth Science I (EarthI),\ Physics IA (PhysIA),\ Earth
		 Science IA (EarthIA)\}
 \end{itemize}
\end{itemize}

Frequencies of the examinees selecting each combination of subjects in 2006
are given in the website of NCUEE. We reproduce part of them in Tables
\ref{tbl:data-soc}--
\ref{tbl:data-sci3} at the end of the paper.
As seen in these tables, examinees may select or not select these
subjects. For example, one examinee may select two subjects from Social 
Studies and three subjects from
Science, while another examinee may select only one subject from 
Science and none from Social Studies. Hence
each examinee is categorized into one of the 
$(6+1)\times\dots\times(4+1)=2800$ 
combinations of individual subjects.
Here 1 is added for not choosing from the subcategory.
As mentioned above, individual departments of universities
impose different additional requirements on the choices of
subjects of NCT.
For example, many science or engineering
departments of national universities
ask the students to take two subjects from Science 
and one subject from Social Studies. 

Let us observe some tendencies of the selections by the examinees 
to illustrate what kind of statistical questions one
might ask concerning the data in Tables
\ref{tbl:data-soc}--\ref{tbl:data-sci3}.
\begin{itemize}
\item[(i)] The most frequent triple of Science subjects is \{BioI,
  ChemI, PhysI\} in Table \ref{tbl:data-sci3}, which
  seems to be consistent with Table \ref{tbl:data-sci} since these three
  subjects are the most frequently selected subjects in 
  Science 1, Science 2 and Science 3,
  respectively. 
  However in Table \ref{tbl:data-sci2}, while the pairs \{BioI,
  ChemI\} and \{ChemI, PhysI\} are the most frequently
  selected pairs in \{Science 1, Science2\} and \{Science 2, Science
  3\}, respectively, the pair \{BioI, PhysI\} is not the
  first choice in \{Science 1, Science 3\}. 
This fact indicates differences 
in the selection of Science subjects between the examinees 
selecting two subjects and those selecting 
three subjects. 

\item[(ii)] In Table
\ref{tbl:data-soc2}  the most frequent
pair is \{GeoB, ContSoc\}. 
However the most frequent single subject 
from Geography and 
History is JHB both in Table \ref{tbl:data-soc} and \ref{tbl:data-soc2}. 
This fact indicates the interaction effect in selecting 
pairs of Social Studies.
\end{itemize}

These observations lead to many  interesting statistical questions.
However Tables \ref{tbl:data-soc}--\ref{tbl:data-sci3} only give
frequencies of choices separately for Social Studies and
Science, i.e., they are the marginal tables for these two major
subjects.
In this paper we are interested in independence across these
two subjects, such as
``are the selections on Social Studies and
Science related or not?''
We give various models for NCT data in Section 3.2 and 
numerical analysis in Section 5.1.

\subsection{The case of Hardy-Weinberg models for allele frequency data}

We also consider problems of population genetics in this paper.
This is another important application of the methodology of this
paper.
The allele frequency data are usually given as the genotype frequency.
For multi-allele locus with alleles $A_1,A_2,\ldots, A_m$, 
the probability of the genotype $A_iA_j$ in an
individual from a random breeding population
is $q_i^2$ $(i = j)$ or $2q_i q_j$ $(i \neq j)$, where
$q_i$ is the proportion of the allele $A_i$.  These are known as the
Hardy-Weinberg equilibrium probabilities. Since the Hardy-Weinberg law
plays an important role in the field of population genetics and often
serves as a basis for genetic inference, much attention has been paid
to tests of the hypothesis that a population being sampled is in the
Hardy-Weinberg equilibrium against the hypothesis that disturbing
forces cause some deviation from the Hardy-Weinberg ratio. See
Crow (1988) and Guo and Thompson (1992)
for example. Though Guo and Thompson (1992) consider the exact test
of the Hardy Weinberg equilibrium for multiple loci, exact procedure
becomes infeasible if the data size or the number of alleles is moderately
large. Therefore MCMC is also useful for this problem. 
Takemura and Aoki (2004) considers conditional tests of 
Hardy-Weinberg model by using MCMC and the technique
of Markov bases.

Due to the rapid progress of sequencing technology, more and more
information is available on the combination of alleles on the same
chromosome.  A combination of alleles at more than one locus on the
same chromosome is called a haplotype and data on haplotype counts are
called haplotype frequency data.  
The haplotype analysis has gained an increasing attention in the
mapping of complex-disease genes, because of the limited power of
conventional single-locus analyses.  
Haplotype data may come with or
without pairing information on homologous chromosomes.
It is technically more difficult to determine pairs of haplotypes
of the corresponding loci on a pair of homologous chromosomes.
A pair of haplotypes on homologous chromosomes is called a diplotype.
In this paper we are interested in diplotype frequency data, because
haplotype frequency data on individual chromosomes without pairing
information are standard contingency table data and can be analyzed by
statistical methods for usual contingency tables.
For the diplotype frequency data,
the null model we want to consider is the independence model that the
probability for each diplotype is expressed by the product of
probabilities for each genotype. 

We consider the models for genotype frequency data in Section 3.3.1
and then consider the models for diplotype
frequency data in Section 3.3.2.
Note that the availability of haplotype data or diplotype data
requires  a separate treatment in our arguments.
Finally we give numerical examples of the analysis of diplotype
frequencies data in Section 5.2.

\section{Conditional tests and models}

\subsection{General formulation of conditional tests and Markov chain
  Monte Carlo procedures}
First we give a brief review on performing MCMC for
conducting conditional tests based on the theory of Markov basis.
Markov basis was introduced by Diaconis and Sturmfels (1998)
and there are now
many references on the definition and the
use of Markov basis (e.g. Aoki and Takemura, 2006).

We denote the space of possible selections as ${\cal I}$.
Each element $\Bi$ in ${\cal I}$ represents a combination of choices.
Following the terminology of contingency tables, each $\Bi \in {\cal I}$
is called  a {\it cell}. 
It should be noted that unlike the
case of standard multiway contingency tables, our index set ${\cal I}$
can not be written as a direct product in general.
We show the structures of ${\cal I}$ for NCT data and allele
frequency data in Section 3.2 and Section 3.3, respectively.

Let $p(\Bi)$ denote the probability of selecting the combination $\Bi$ (or
the probability of cell $\Bi$) and write $\Bp = \{p(\Bi)\}_{\Bi \in
  {\cal I}}$.  
In this paper, we do not necessarily assume that $\Bp$ is
normalized. In fact, in the models we consider in this paper, we only 
give an unnormalized functional specification of $p(\cdot)$.
Note that we need not calculate the normalizing 
constant $\sum_{\Bi\in {\cal I}} p(\Bi)$ for performing a MCMC procedure.
Denote the result of the selections by $\samplesize$ individuals 
as $\Bx = \{x(\Bi)\}_{\Bi \in {\cal I}}$,
where $x(\Bi)$ is the frequency of the cell $\Bi$.  We call $\Bx$ a
frequency vector. 

In the models considered in this paper, the cell probability $p(\Bi)$
is written as some product of functions, which correspond to various
marginal probabilities. 
Let ${\cal J}$ denote the index  set of the marginals.
Then our models can be written as
\begin{equation}
\label{eq:general-model}
p(\Bi)= h(\Bi) \prod_{\Bj\in {\cal J}} q(\Bj)^{a_{\Bj \Bi}},
\end{equation}
where $h(\Bi)$ is a known function and $q(\Bj)$'s are the parameters.
An important point here is that 
the sufficient statistic $\Bt=\{t(\Bj), \Bj\in {\cal J}\}$
is written in a matrix form as
\begin{equation}
\Bt = A \Bx,\quad A=(a_{\Bj\Bi})_{\Bj\in{\cal J}, \Bi\in{\cal I}},
\label{eq:t=Ax}
\end{equation}
where $A$ is $\numrowsA \times \numcolsA$ matrix of non-negative
integers and $\numrowsA=|{\cal J}|$,  $\numcolsA=|{\cal I}|$.  
We call $A$ a {\em
  configuration} in connection with the theory of toric ideals
in Section 4.

By the standard theory of conditional tests (Lehmann and Romano, 2005,
for example), we can
perform conditional test of the model (\ref{eq:general-model}) based
on the conditional distribution given the sufficient statistic $\Bt$.
The conditional sample space given $\Bt$,  called the $\Bt$-fiber, 
is 
\[
{\cal F}_{\Bt} = \{ \Bx \in \NN^\numcolsA \mid \Bt = A \Bx\},
\]
where $\NN=\{0,1,\dots\}$.
If we can sample from the 
conditional distribution over ${\cal F}_{\Bt}$, we can evaluate
$P$-values of any test statistic.  One of the advantages of 
MCMC method of sampling is that it can be run without evaluating the
normalizing constant.
Also once a connected Markov chain over the conditional
sample space is constructed, then 
the chain can be modified to
give a connected and aperiodic Markov chain with the stationary
distribution 
by the Metropolis-Hastings
procedure (e.g. Hastings, 1970). 
Therefore it is essential to construct a connected
chain and the solution to this problem is given by the notion of
{\em Markov basis} (Diaconis and Sturmfels, 1998).

The fundamental contribution of 
Diaconis and Sturmfels (1998) 
is to show that a Markov basis is given as 
a binomial generator of the
well-specified polynomial ideal (toric ideal) and it can be 
given as a Gr\"obner basis. 
In Section 4, we show that our problem 
considered in Section 3.2 and 3.3
corresponds to a 
well-known toric ideal and give an explicit form of the 
reduced Gr\"obner basis.

\subsection{Models for NCT data}
\label{sec:NCT-models}

Following the general formalization in Section 3.1, we 
formulate data types and their statistical models
in view  of NCT.  
Suppose that there are $\numgroup$ different groups (or categories)
and  $\numsubgroup_j$ different 
subgroups in group $j$ for  $j = 1,\ldots,\numgroup$.
There are $\numjkn_{jk}$ different {\it items} in subgroup $k$ of group $j$
($k = 1,\ldots,\numsubgroup_j$, $j = 1,\ldots,\numgroup$).
In NCT, $\numgroup=2$,
$\numsubgroup_1 = |\{\mbox{Geography and History, Civics}\}| = 2$ 
and similarly  $\numsubgroup_2=3$.
The sizes of subgroups are
$\numjkn_{11} =  | \{ \text{WHA, WHB, JHA, JHB},\allowbreak \text{GeoA,
  GeoB}\} | = 6$ and
similarly 
$\numjkn_{12} = 3$,
$\numjkn_{21} = 4$,
$\numjkn_{22} = 3$,
$\numjkn_{23}=4$.

Each individual selects $\numjkc_{jk}$ items from 
the subgroup $k$ of group $j$.  We assume that the total number
$\totalitems$ 
of items  chosen is fixed and common for all individuals.  
In NCT $\numjkc_{jk}$ is either 0 or 1.  
For example if an examinee is required to take two Science subjects in 
NCT, then $(\numjkc_{21},\numjkc_{22},\numjkc_{23})$ is
$(1,1,0)$, $(1,0,1)$ or $(0,1,1)$.
For the analysis of genotypes  in Section \ref{sec:HW},
$\numjkc_{jk}\equiv 2$ although there is no nesting of
subgroups, and  the same item (allele) can be selected more than once
(selection ``with replacement'').

We now set up our notation for indexing a combination of choices 
somewhat carefully.  In NCT, if an examinee chooses 
WHA from ``Geography and History'' of Social Studies and PhysI from
Science 3 of Science, we denote the combination of these two choices as
(111)(231).   In this notation, 
the selection of $\numjkc_{jk}$ items from the subgroup $k$ of group $j$
are indexed as
\[
\Bi_{jk}=(jkl_1)(jkl_2)\dots(jkl_{\numjkc_{jk}}), \qquad
1 \le l_1 \le \dots \le l_{\numjkc_{jk}} \le \numjkn_{jk}.
\]
Here $\Bi_{jk}$ is regarded as a string.
If nothing is selected from the subgroup, we define 
$\Bi_{jk}$ to be an empty string.  
Now  by concatenation of strings, 
the set $\cal I$ of combinations is written as
\[
{\cal I} = \{\Bi=\Bi_1\dots  \Bi_\numgroup \},
\qquad \Bi_j = \Bi_{j1}\dots \Bi_{j\numsubgroup_j}, 
\quad  j=1,\dots,\numgroup. 
\]
For example the choice of (P\&E, BioI, ChemI) in NCT is denoted by
$\Bi=(123)(212)(222)$.  
In the following we denote $\Bi'\subset \Bi$ if $\Bi'$ appears as a
substring of $\Bi$.

Now we consider some statistical models for $\Bp$.
For NCT data, we consider three simple statistical models, namely, {\it
complete 
independence model}, {\it subgroup-wise independence model}
and {\it group-wise independence model}.
The complete independence model is defined as
\begin{equation}
\label{eq:complete-independece}
p(\Bi) 
= \displaystyle\prod_{j = 1}^\numgroup
  \prod_{k = 1 \atop \Bi_{jk} \subset \Bi}^{\numsubgroup_j} 
\prod_{t = 1}^{\numjkc_{jk}}
\marginalprobs_{jk}(l_t)
\end{equation}
for some parameters $\marginalprobs_{jk}(l),\ j = 1,\ldots,\numgroup;\ k =
1,\ldots,m_j;\ l= 1,\ldots,\numjkn_{jk}$.
Note that if $\numjkc_{jk} >1$ we need a multinomial coefficient in
(\ref{eq:complete-independece}).
The complete independence 
model means that each
$p(\Bi)$, the inclination of the combination $\Bi$, 
is explained by the set of 
inclinations $q_{jk}(l)$ of each item.
Here $\marginalprobs_{jk}(l)$ corresponds to the marginal
probability of the item $(jkl)$.  However we  do not necessarily normalize them
as $1=\sum_{l=1}^{\numjkn_{jk}}\marginalprobs_{jk}(l)$, because the
normalization for $\Bp$ is not trivial anyway.  The same comment
applies to other models below.

Similarly, 
the subgroup-wise independence model is defined as
\begin{equation}
\label{eq:subgroup-indepedence}
p(\Bi) 
 = \displaystyle\prod_{j = 1}^\numgroup\prod_{k = 1 \atop
 \Bi_{jk} \subset \Bi }^{\numsubgroup_j}
\marginalprobs_{jk}(\Bi_{jk})
\end{equation}
for some parameters $\marginalprobs_{jk}(\cdot)$,
and the group-wise independence model is defined as
\begin{equation}
\label{eq:group-independence}
p(\Bi) 
 = \displaystyle\prod_{j = 1}^\numgroup
\marginalprobs_j(\Bi_j)
\end{equation}
for some parameters $\marginalprobs_j(\cdot)$.

In this paper, we treat these models as the {\it null models} and give
testing procedures 
to assess their fitting to observed data following the general theory in
Section 3.1.

\subsection{Models for allele frequency data}
\label{sec:HW}

\subsubsection{Models for the genotype frequency data}
\label{sec:genotype}
We assume that there are $\numgroup$ distinct loci.
In the locus $j$, there are $m_j$ distinct alleles, 
$A_{j1},\ldots,A_{jm_j}$.
In this case, we can imagine that
each individual selects two alleles for 
each locus {\it with replacement}. 
Therefore the set of the combinations is  written as
\[
\begin{array}{c}
{\cal I} = \{\Bi =  (i_{11}i_{12})(i_{21}i_{22})\ldots
(i_{\numgroup 1}i_{\numgroup 2})
\  |\ 1\leq i_{j1} \leq i_{j2} \leq m_{j},\ j = 1,\ldots,\numgroup\}.
\end{array}
\]

For the genotype frequency data, we consider two models of
hierarchical structure, namely,
{\it genotype-wise independence model}
\begin{equation}
\label{eq:genotype-indepedent}
p(\Bi) = 
\prod_{j = 1}^\numgroup \marginalprobs_j(i_{j1}i_{j2})
\end{equation}
and the Hardy-Weinberg model
\begin{equation}
\label{eq:genetype-HW}
p(\Bi) = \prod_{j = 1}^\numgroup \tilde \marginalprobs_j(i_{j1}i_{j2}),
\end{equation}
where
\begin{equation}
\label{eq:HW2}
\tilde \marginalprobs_j(i_{j1}i_{j2}) = \left\{\begin{array}{ll}
q_j(i_{j1})^2 & \mbox{if}\ i_{j1} = i_{j2},\\
2q_j(i_{j1})q_j(i_{j2}) & \mbox{if}\ i_{j1} \neq i_{j2}.
\end{array}\right.
\end{equation}
Note that for both cases the sufficient statistic $\Bt$ can be written
as $\Bt = A \Bx$ for an appropriate matrix $A$ as shown in Section
\ref{subsec:diplo-data}.

\subsubsection{Models for the diplotype frequency data}
\label{sec:haplotype}
In order to illustrate the difference between genotype data and
diplotype data, consider a simple case of $\numgroup=2, \numsubgroup_1
=\numsubgroup_2 = 2$ and
suppose that genotypes of $\samplesize  = 4$ individuals are given as
\[
\{A_{11}A_{11},A_{21}A_{21}\}, \ 
\{A_{11}A_{11},A_{21}A_{22}\}, \ 
 \{A_{11}A_{12},A_{21}A_{21}\}, \ 
\{A_{11}A_{12},A_{21}A_{22}\}.
\]
In this genotype data, for an individual who has homozygote genotype
on at least one loci, the diplotypes are uniquely determined. However,
for the fourth individual who has the genotype
$\{A_{11}A_{12},A_{21}A_{22}\}$, there are two possible 
diplotypes as $\{(A_{11},A_{21}), (A_{12},A_{22})\}$ and
$\{(A_{11},A_{22}), (A_{12},A_{21})\}$.

Now suppose that  information on diplotypes are available.
The set of combinations for the
diplotype data is given as
\[
\begin{array}{c}
{\cal I} = \{\Bi =  \Bi_1 \Bi_2=
(i_{11}\cdots i_{\numgroup 1})(i_{12}\cdots i_{\numgroup 2})
\  |\ 1\leq i_{j1} , i_{j2} \leq m_{j},\ j = 1,\ldots,\numgroup\}.
\end{array}
\]
In order to determine the order of $\Bi_1=(i_{11}\dots i_{r1})$ and
$\Bi_2=(i_{12}\dots i_{r2})$ uniquely, we assume that these two are
lexicographically ordered, i.e., there exists some $j$ 
such that 
\[
i_{11}=i_{12}, \dots, i_{j-1,1}=i_{j-1,2}, \ i_{j1}< i_{j2}
\]
unless $\Bi_1 = \Bi_2$.

For the parameter $\Bp = \{p(\Bi)\}$ where $p(\Bi)$ is the probability
for the diplotype $\Bi$, we can consider the same models as for the
genotype case.  Corresponding to the null hypothesis that
diplotype data do not contain more information than the genotype data, 
we can consider the genotype-wise independence model 
(\ref{eq:genotype-indepedent}) 
and the Hardy-Weinberg model 
(\ref{eq:genetype-HW}).  The sufficient statistics for these models
are the same as in the previous subsection.

If these models are rejected, we can further
test independence in diplotype data.  For example we can consider
a haplotype-wise Hardy-Weinberg model.
\[
p(\Bi) = p(\Bi_1 \Bi_2)=
\left\{\begin{array}{ll}
q(\Bi_1)^2 & \mbox{if}\ \Bi_1 = \Bi_2,\\
2 q(\Bi_1) q(\Bi_2) & \mbox{if}\ \Bi_1 \neq \Bi_2.
\end{array}\right.
\]
The sufficient statistic for this model is given by the set 
of frequencies of each haplotype and the conditional test can be
performed as in the case of Hardy-Weinberg model for a single gene
by formally identifying each haplotype as an allele.

\section{Gr\"obner basis for Segre-Veronese configuration}
\label{sec:GB}

In this section,
we introduce toric ideals of algebras of Segre-Veronese type
(Ohsugi and Hibi, 2000) with a generalization 
to fit statistical applications in the present paper.

First we define toric ideals.
A {\em configuration} in $\RR^\numrowsA$ is a finite set 
$\Ac =\{\ab_1,\ldots,\ab_\numcolsA \} \subset \NN^d$.
$A$ can be regarded as a $\numrowsA \times \numcolsA$ matrix and
corresponds to the matrix connecting the frequency vector to 
the sufficient statistic as in (\ref{eq:t=Ax}).
Let $K$ be a field and $K[\qb]=K[q_1, \dots, q_\numrowsA]$
the polynomial ring in $\numrowsA$
variables over $K$.
We associate a configuration $\Ac \subset \NN^d$
with the semigroup ring
$K[\Ac] = K[\qb^{\ab_1}, \ldots, \qb^{\ab_\numcolsA}]$
where $\qb^{\ab} = q_1^{a_1} \cdots q_d^{a_d}$
if $\ab = (a_1, \ldots, a_d)$.
Note that $\numrowsA=|{\cal J}|$ and $\qb^{\ab_i}$ corresponds to
to the term $\prod_{\Bj\in{\cal J}}q(\Bj)^{a_{\Bj \Bi}}$ on the
right-hand side of 
(\ref{eq:general-model}).
Let $K[W] = K[w_1,\ldots,w_\numcolsA]$ be
the polynomial ring in $\numcolsA$ variables over $K$.
Here $\numcolsA=|{\cal I}|$ and the variables $w_1,\ldots,w_\numcolsA$
correspond to the cells of ${\cal I}$.
The {\em toric ideal} $I_\Ac$ of $\Ac$ 
is the kernel of the surjective homomorphism 
$\pi \, : \, K[W] \to K[\Ac]$
defined by setting
$\pi(w_i) = \qb^{\ab_i}$ for all $1 \leq i \leq \numcolsA$.
It is known that the toric ideal
$I_\Ac$ is generated by the binomials $u - v$, where
$u$ and $v$ are monomials of $K[W]$, with
$\pi(u) = \pi(v)$.
More precisely, $I_\Ac$ is written as
$$
I_\Ac = 
\left<
\left.
W^{\zb^{+}} - 
W^{\zb^{-}} \ \right| \ \zb \in \ZZ^\numcolsA, \ 
\Ac \zb = {\bf 0}
\right>,
$$
where $\zb = \zb^{+} -\zb^{-}$
with $\zb^{+}, \zb^{-} \in \NN^\numcolsA$.
We call an integer vector $\Bz\in \ZZ^\numcolsA$ a {\it move} if
$\Ac \zb ={\bf 0}$.

The  {\it initial ideal} of $I_A$
with respect to a monomial order is
the ideal of $K[W]$
generated by all initial monomials of nonzero elements of $I_A$.
A finite set ${\cal G}$ of $I_A$ is called a
{\it Gr\"obner basis}
of $I_A$ with respect to a monomial order $<$
if the initial ideal of $I_A$ with respect to
$<$ is generated by the initial monomials of 
the polynomials in ${\cal G}$.
A Gr\"obner basis ${\cal G}$ is called {\it reduced}
if, for each $g \in {\cal G}$, none of the monomials in $g$ is divisible by
the initial monomials of $g'$ for some $g \neq g' \in {\cal G}$.
It is known that
if ${\cal G}$ is a Gr\"obner basis of $I_A$,
then $I_A$ is generated by ${\cal G}$.
In general, the reduced Gr\"obner basis of a toric ideal
consists of binomials.
See Chapter 4 of Sturmfels (1995) for the details of toric ideals and Gr\"obner bases.

The following proposition
associates
Markov bases with toric ideals.

\begin{proposition}[Diaconis--Sturmfels, 1998]
A set of moves ${\cal B} =
 \{\zb_1,\ldots,\zb_L\}$
is a Markov basis if and only if
$I_A$ is generated by binomials
$W^{\zb_1^+} - W^{\zb_1^-}$,
$\ldots$,
$W^{\zb_L^+} - W^{\zb_L^-}$.
\end{proposition}

We now introduce the notion of 
algebras of Segre-Veronese type.
Fix integers $\totalitems \geq 2$, $M\ge 1$ and
sets of integers
${\bf b} = \{b_1,\ldots,b_M\}$,
${\bf c} = \{c_1,\ldots,c_M\}$,
${\bf r} = \{r_1,\ldots,r_M\}$ and
${\bf s} = \{s_1,\ldots,s_M\}$ such that 
\begin{enumerate}
\item[(i)]
$0 \leq c_i \leq b_i$
for all $1 \leq i \leq M$;
\item[(ii)]
$1 \leq s_i \leq r_i \leq \numrowsA$
for all $1 \leq i \leq M$.
\end{enumerate}
Let $\AAA \subset \NN^{\numrowsA}$ denote the configuration
consisting of all nonnegative integer vectors
$(f_1,f_2,\ldots,f_\numrowsA)
\in \NN^{\numrowsA}$
such that
\begin{enumerate}
\item[(i)]
$\sum_{j=1}^{d} f_j =\totalitems$.
\item[(ii)]
$c_i \leq \sum_{j=s_i}^{r_i} f_j \leq  b_i$
for all $1 \leq i \leq M$.
\end{enumerate}
Let $K[\AAA]$ denote the affine semigroup ring 
generated by all monomials 
$\prod_{j=1}^{\numrowsA} {q_j}^{f_j}$
over $K$
and call it an
{\it algebra of Segre-Veronese type}.
Note that the present definition generalizes
the definition in Ohsugi and Hibi (2000).

Several
popular classes of semigroup rings are
algebras of Segre-Veronese type.
If $M=2$, $\totalitems=2$, $b_1= b_2=c_1=c_2=1$,
$s_1 = 1$, $s_2 = r_1+1$ and $r_2 =\numrowsA$, then
the affine semigroup ring
$K[\AAA]$
is
the Segre product of polynomial rings
$K[q_1,\ldots,q_{r_1}]$
and
$K[q_{r_1+1},\ldots,q_\numrowsA]$.
On the other hand,
if $M=\numrowsA$, $s_i=r_i=i$, $b_i=\totalitems$ and $c_i = 0$
for all $1 \leq i \leq M$, then the affine semigroup ring
$K[\AAA]$
is the classical
$\totalitems$th Veronese subring of the polynomial ring
$K[q_1,\dots,q_\numrowsA]$.
Moreover, 
if $M=\numrowsA$, $s_i=r_i=i$, $b_i=1$ and $c_i = 0$
for all $1 \leq i \leq M$, then the affine semigroup ring
$K[\AAA]$
is the $\totalitems$th 
squarefree Veronese subring of the polynomial ring
$K[q_1,\dots,q_\numrowsA]$.
In addition, algebras of Veronese type 
(i.e., $M=\numrowsA$, $s_i=r_i=i$ and $c_i = 0$
for all $1 \leq i \leq M$)
are 
studied in De Negri and Hibi (1997) and Sturmfels (1995).

Let $K[Y]$ denote the polynomial ring with 
the set of variables 
$$
\left\{ y_{j_1 j_2 \cdots j_\totalitems}
\ \left| \  1\leq j_1 \leq  j_2 \leq \cdots \leq j_\totalitems
\leq d, \ \prod_{k=1}^\totalitems q_{j_k} 
\in \{ \qb^{\ab_1} , \ldots, \qb^{\ab_\numcolsA}\}
\right. \right\}, 
$$
where $K[\AAA]=K[\qb^{\ab_1}, \ldots, \qb^{\ab_\numcolsA}]$. 
The toric ideal $I_{\AAA}$
is the kernel of the surjective homomorphism
$
\pi : K[Y]
\longrightarrow 
K[\AAA]
$
defined by
$
\pi(y_{j_1 j_2 \cdots j_\totalitems})
=
\prod_{k=1}^\totalitems q_{j_k}
$.

A monomial
$
y_{\alpha_1 \alpha_2 \cdots \alpha_\totalitems}
y_{\beta_1 \beta_2 \cdots \beta_\totalitems}
\cdots
y_{\gamma_1 \gamma_2 \cdots \gamma_\totalitems}
$
is called {\it sorted} if 
$$
\alpha_1 \leq \beta_1 \leq \cdots \leq \gamma_1 \leq 
\alpha_2 \leq \beta_2 \leq \cdots \leq \gamma_2 \leq \cdots \leq 
\alpha_\totalitems \leq \beta_\totalitems \leq \cdots \leq \gamma_\totalitems
.$$
Let ${\rm sort}(\cdot)$ denote
the operator which takes any string over
the alphabet $\{1,2,\ldots,d\}$
and sorts it into weakly increasing order.
Then the quadratic Gr\"obner basis of toric
ideal $I_\AAA$ is given as follows.

\begin{theorem}
\label{thm:GB}
Work with the same notation as above.
Then there exists a monomial order on $K[Y]$
such that the set of all binomials 
\begin{equation}
\label{eq:binomial}
\{
y_{\alpha_1 \alpha_2 \cdots \alpha_\totalitems}
y_{\beta_1 \beta_2 \cdots \beta_\totalitems}
-
y_{\gamma_1 \gamma_3 \cdots \gamma_{2\totalitems-1}}
y_{\gamma_2 \gamma_4 \cdots \gamma_{2\totalitems}}
\ 
|
\ 
{\rm sort}(
{\alpha_1 \beta_1 \alpha_2 \beta_2 \cdots \alpha_\totalitems \beta_\totalitems}
)
=
\gamma_1 \gamma_2 \cdots \gamma_{2\totalitems}
\}
\end{equation}
is the reduced Gr\"obner basis of
the toric ideal $I_\AAA$.
The initial ideal is generated by 
squarefree quadratic (nonsorted) monomials.

In particular, 
the set of all integer vectors corresponding to
the above binomials is a Markov basis.  Furthermore the set is minimal as a 
Markov basis.
\end{theorem}

\noindent
{\it Proof.}
The basic idea of the proof appears in Theorem 14.2 in Sturmfels (1995).

Let ${\cal G}$ be the above set of binomials.
First we show that ${\cal G} \subset I_\AAA$.
Suppose that $m = y_{\alpha_1 \alpha_2 \cdots \alpha_\totalitems}
y_{\beta_1 \beta_2 \cdots \beta_\totalitems}
$ is not sorted and let
$$
\gamma_1 \gamma_2 \cdots \gamma_{2\totalitems}
=
{\rm sort}(
{\alpha_1 \beta_1 \alpha_2 \beta_2 \cdots \alpha_\totalitems \beta_\totalitems}
)
.$$
Then,
$m$ is squarefree
since
the monomial
$y_{\alpha_1 \alpha_2 \cdots \alpha_\totalitems}^2$
is sorted.
Since the binomial
$
y_{\alpha_1 \alpha_2 \cdots \alpha_\totalitems}
y_{\beta_1 \beta_2 \cdots \beta_\totalitems}
-
y_{\alpha_1' \alpha_2' \cdots \alpha_\totalitems'}
y_{\beta_1' \beta_2' \cdots \beta_\totalitems'}
\in K[Y]$
belongs to $I_{\AAA}$ if and only if
$
{\rm sort}(\alpha_1 \alpha_2 \cdots \alpha_\totalitems \beta_1 \beta_2 \cdots \beta_\totalitems)
=
{\rm sort}(\alpha_1' \alpha_2' \cdots \alpha_\totalitems' \beta_1' \beta_2' \cdots \beta_\totalitems')
$,
it is sufficient to show that both
$y_{\gamma_1 \gamma_3 \cdots \gamma_{2\totalitems-1}}$
and
$y_{\gamma_2 \gamma_4 \cdots \gamma_{2\totalitems}}$
are variables of $K[Y]$.
For $1 \leq i \leq n$, let
$
\rho_i  =  | \{ j \ | \ s_i \leq \gamma_{2 j -1} \leq r_i\} |$
and
$\sigma_i =   | \{ j \ | \ s_i \leq \gamma_{2 j  } \leq r_i\} |
$. Since 
$\gamma_1 \leq \gamma_2 \leq \cdots \leq \gamma_{2\totalitems}$,
$\rho_i$ and $\sigma_i$
are either equal or they differ by one
for each $i$.
If $\rho_i \leq \sigma_i$,
then $0 \leq \sigma_i -\rho_i \leq 1$.
Since $2 c_i \leq \rho_i + \sigma_i \leq 2 b_i$,
we have $\sigma_i \leq b_i + 1/2$
and $c_i - 1/2 \leq \rho_i$.
Thus $c_i \leq \rho_i \leq \sigma_i \leq b_i$.
If 
$\rho_i > \sigma_i$,
then $\rho_i -\sigma_i = 1$.
Since $2 c_i \leq \rho_i + \sigma_i \leq 2 b_i$,
we have $\rho_i \leq b_i + 1/2$
and $c_i - 1/2 \leq \sigma_i$.
Thus $c_i \leq \sigma_i < \rho_i \leq b_i$.
Hence 
$y_{\gamma_1 \gamma_3 \cdots \gamma_{2\totalitems-1}}$
and
$y_{\gamma_2 \gamma_4 \cdots \gamma_{2\totalitems}}$
are variables of $K[Y]$.

By virtue of relation
between
the reduction of a monomial by ${\cal G}$ and 
sorting of the indices of a monomial,
it follows that
there exists a monomial order such that,
for any binomial in ${\cal G}$,
the first monomial is the initial monomial.
See also Theorem 3.12 in Sturmfels (1995).

Suppose that ${\cal G}$ is not a Gr\"obner basis.
Thanks to Macaulay's Theorem,
there exists a binomial $f \in I_{\AAA}$ such that
both monomials in $f$ are sorted.
This means that $f = 0$
and $f$ is not a binomial. 
Hence ${\cal G}$ is a Gr\"obner basis of $I_{\AAA}$.
It is easy to see that the Gr\"obner basis ${\cal G}$ is reduced
and a minimal set of generators of $I_\AAA$.
\hfill{Q.E.D.}

\bigskip

Finally we describe how to run a Markov chain using the Gr\"obner basis
given in Theorem \ref{thm:GB}. First, given a configuration $A$ in
(\ref{eq:t=Ax}), we check that (with appropriate reordering of rows)
that $A$ is indeed a configuration of Segre-Veronese type.  It is easy
to check that our models in Sections \ref{sec:NCT-models} and
\ref{sec:HW} are of Segre-Veronese type, because the restrictions on
choices are imposed separately for each group or each subgroup. Recall
that each column of $A$ consists of non-negative integers  whose sum
$\totalitems$ is common.

We now associate to each column $\ab_i$ of $A$ a
set of indices indicating the rows with  positive elements $a_{\Bj
  \Bi}>0$ and a particular index $\Bj$ is repeated $a_{\Bj \Bi}$
times.  For example if $\numrowsA=4, \totalitems=3$ and
$\ab_i=(1,0,2,0)'$, then row 1 appears once and 
row 3 appears twice in $\ab_i$.  Therefore we associate the index
$(1,3,3)$ to $\ab_i$.  We can consider the set of indices as
$\totalitems \times \numcolsA$ matrix $\tilde A$.  Note that $\tilde
A$ and $A$ carry the same information.

Given $\tilde A$, we can choose a random element of the reduced
Gr\"obner basis of Theorem \ref{thm:GB} as follows.  Choose two
columns (i.e.\ choose two cells from $\cal I$)
of $\tilde A$ and sort $2\times\totalitems$ elements of these
two columns.  {}From the sorted elements, pick alternate elements and
form two new sets of indices.   For example
if $\totalitems=3$ and the two chosen columns of $\tilde A$ are 
$(1,3,3)$ and $(1,2,4)$, then by sorting these 6 elements we obtain 
$(1,1,2,3,3,4)$.  Picking alternate elements produces
$(1,2,3)$ and $(1,3,4)$.  These new sets of indices correspond to (a
possibly overlapping) two columns of $\tilde A$, hence to two cells of
$\cal I$.  Now the difference of the two original columns and the two
sorted columns of $\tilde A$ correspond to a random binomial  in
(\ref{eq:binomial}). 
It should be noted that when the sorted columns coincide with the
original columns, then we discard these columns and choose other two
columns. The rest of
the procedure for running a Markov chain is described in Diaconis and Sturmfels (1998).
See also Aoki and Takemura (2006).

\section{Numerical examples}
In this section we present numerical experiments on NCT data and
a diplotype frequency data.

\subsection{The analysis of NCT data}

First we consider the analysis of NCT data concerning selections in
Social Studies and Science.  
Because NCUEE currently do not provide cross tabulations of
frequencies of choices across the major subjects, we can not evaluate
the $P$-value of the actual data.  However for the models
in Section \ref{sec:NCT-models}, the sufficient statistics (the
marginal frequencies) can be obtained from Tables
$\ref{tbl:data-soc}$--
\ref{tbl:data-sci3}.  Therefore in this section we evaluate the
conditional null distribution of the Pearson's $\chi^2$  statistic by MCMC
and compare it to the asymptotic $\chi^2$ distribution.

In Section \ref{sec:NCT-models}, we consider 
three models, complete independence model, subgroup-wise independence model and
group-wise independence model, for the setting of group-wise selection
problems.
Note that, however, the subgroup-wise independence model coincides with
the group-wise independence model
for NCT data, since $c_{jk} \leq 1$ for all $j$ and $k$. Therefore we
consider fitting of the complete
independence model and the group-wise independence model for NCT data.

As we have seen in Section 2.1, there are many kinds of
choices for each examinee.
However, it may be natural to treat some similar subjects as one
subject. For example, WHA and WHB
may well be treated as WH, ChemI and Chem IA may well be treated as
Chem, and so on.
As a result, we
consider the following aggregation of subjects.
\begin{itemize}
\item In Social Studies:
WH = \{WHA,WHB\}, JH = \{JHA,JHB\}, Geo = \{GeoA,GeoB\}
\item In Science:
CSiB = \{CSiB, ISci\},  Bio = \{BioI, BioIA\}, 
Chem = \{ChemI, ChemIA\}, Phys = \{PhysI, PhysIA\},
Earth = \{EarthI, EarthIA\}
\end{itemize}
In our analysis, we take a look at examinees selecting 
two subjects for Social Studies and two subjects for Science.
Therefore 
\[\begin{array}{l}
J = 2, m_1 = 2, m_2 = 3, m_{11} = m_{12} = 3, m_{21} = m_{22} = m_{23} =
 2,\\
c_{11} = c_{12} = 1, (c_{21},c_{22},c_{23}) = (1,1,0)\ \mbox{or}\
 (1,0,1)\ \mbox{or}\ (0,1,1).
\end{array}
\]
The number of possible combination is then
$\nu = |{\cal I}| = 3\cdot 3\times 3\cdot2^2 = 108$.
Accordingly our sample size is $n=195094$, 
which is the number of examinees selecting two subjects
on Science from Table \ref{tbl:data-sci}.  Our data set is shown
in Table~\ref{tbl:data-NCT-imaginary}.

\begin{table*}[tbhp]
\begin{center}
\caption{The data set of number of the examinees in NCT in 2006 ($n = 195094$)}
\label{tbl:data-NCT-imaginary}
\begin{tabular}{|c|c|c|c|} \hline
\multicolumn{1}{|c}{} & \multicolumn{1}{|c|}{ContS} &
 \multicolumn{1}{c|}{Ethics} & \multicolumn{1}{c|}{P\&E}\\ \hline
WH & 32352 & 8839 & 8338\\ \hline
JH & 51573 & 8684 & 14499\\ \hline
Geo & 59588 & 4046 & 7175\\ \hline
\end{tabular}
\hspace*{1cm}
\begin{tabular}{|c|c|c|c|c|}\hline
& CSiA & Chem & Phys & Earth\\ \hline
CSiB  & 1648 & 1572 & 169 & 4012\\ \hline
Bio   & 21392 & 55583 & 1416 & 1845\\ \hline
Phys  & 3286 & 102856 & --- & ---\\ \hline
Earth & 522 & 793 & --- & ---\\ \hline
\end{tabular}
\end{center}
\end{table*}

From Table \ref{tbl:data-NCT-imaginary}, we can calculate the maximum
likelihood estimates
of the numbers of the examinees selecting each combination of subjects. 
The sufficient statistics under the complete independence model
are the numbers of the examinees selecting each subject, whereas
the sufficient statistics under the group-wise independence model
are the numbers of the examinees selecting each combination of subjects
in the same group.
The maximum likelihood estimates calculated from the sufficient
statistics are shown in 
Table \ref{tbl:MLE-NCT}.  For the complete independence model the
maximum likelihood estimates  can
be calculated as in Section 5.2 of Bishop {\it et al.} (1975).

\begin{table*}[htbp]
\begin{center}
\caption{MLE of the number of the examinees selecting each combination of
subjects under the complete independence model (upper) and the group-wise independence model (lower).}
\label{tbl:MLE-NCT}
{\tiny
\begin{tabular}{|c|c|c|c|c|c|c|c|c|c|}\hline
 & \multicolumn{3}{|c|}{WH} & \multicolumn{3}{|c|}{JH} & \multicolumn{3}{|c|}{Geo}\\ \cline{2-10}
 & ContS & Ethics & P\&E & ContS & Ethics & P\&E & ContS & Ethics & P\&E \\ \hline
CSiB,CSiA  &  180.96 &    27.20 &    37.84 &   273.12 &    41.05 &    57.12 &   258.70 &    38.88 &    54.10\\ 
           &  273.28 &    74.66 &    70.43 &   435.65 &    73.36 &   122.48 &   503.35 &    34.18 &    60.61\\ \hline
CSiB,Chem  & 1083.82 &   162.89 &   226.65 &  1635.85 &   245.86 &   342.10 &  1549.48 &   232.88 &   324.03\\ 
           &  260.68 &    71.22 &    67.18 &   415.56 &    69.97 &   116.83 &   480.14 &    32.60 &    57.81\\ \hline
CSiB,Phys  &  110.04 &    16.54 &    23.01 &   166.09 &    24.96 &    34.73 &   157.32 &    23.64 &    32.90\\ 
           &   28.02 &     7.66 &     7.22 &    44.68 &     7.52 &    12.56 &    51.62 &     3.50 &     6.22\\ \hline
CSiB,Earth &    7.33 &     1.10 &     1.53 &    11.06 &     1.66 &     2.31 &    10.47 &     1.57 &     2.19\\ 
           &  665.30 &   181.77 &   171.47 &  1060.57 &   178.58 &   298.16 &  1225.39 &    83.20 &   147.55\\ \hline
Bio,CSiA   & 1961.78 &   294.84 &   410.26 &  2960.99 &   445.02 &   619.21 &  2804.66 &   421.52 &   586.52\\ 
           & 3547.39 &   969.19 &   914.26 &  5654.96 &   952.20 &  1589.81 &  6533.81 &   443.64 &   786.74\\ \hline
Bio,Chem   &11749.94 &  1765.93 &  2457.19 & 17734.63 &  2665.39 &  3708.74 & 16798.27 &  2524.66 &  3512.92\\ 
           & 9217.20 &  2518.26 &  2375.53 & 14693.34 &  2474.10 &  4130.82 & 16976.84 &  1152.72 &  2044.18\\ \hline
Bio,Phys   & 1193.01 &   179.30 &   249.49 &  1800.65 &   270.63 &   376.56 &  1705.58 &   256.34 &   356.68\\ 
           &  234.81 &    64.15 &    60.52 &   374.32 &    63.03 &   105.23 &   432.49 &    29.37 &    52.08\\ \hline
Bio,Earth  &   79.43 &    11.94 &    16.61 &   119.88 &    18.02 &    25.07 &   113.55 &    17.07 &    23.75\\ 
           &  305.95 &    83.59 &    78.85 &   487.72 &    82.12 &   137.12 &   563.52 &    38.26 &    67.85\\ \hline
CSiA,Phys  & 2691.94 &   404.58 &   562.95 &  4063.04 &   610.65 &   849.68 &  3848.52 &   578.41 &   804.82\\ 
           &  544.91 &   148.88 &   140.44 &   868.65 &   146.27 &   244.21 &  1003.65 &    68.15 &   120.85\\ \hline
CSiA,Earth &  179.22 &    26.94 &    37.48 &   270.50 &    40.65 &    56.57 &   256.22 &    38.51 &    53.58\\ 
           &   86.56 &    23.65 &    22.31 &   137.99 &    23.24 &    38.79 &   159.44 &    10.83 &    19.20\\ \hline
Bio,Phys   &16123.14 &  2423.20 &  3371.73 & 24335.27 &  3657.42 &  5089.09 & 23050.40 &  3464.31 &  4820.39\\ 
           &17056.38 &  4660.03 &  4395.90 & 27189.93 &  4578.31 &  7644.05 & 31415.54 &  2133.10 &  3782.75\\ \hline
Bio,Earth  & 1073.41 &   161.33 &   224.48 &  1620.14 &   243.50 &   338.81 &  1534.60 &   230.64 &   320.92\\ 
           &  131.50 &    35.93 &    33.89 &   209.63 &    35.30 &    58.93 &   242.21 &    16.45 &    29.16\\ \hline
\end{tabular}
}
\end{center}
\end{table*}

The configuration $A$ for the
complete independence model is written as
\[
A = \left[
\begin{array}{ccc}
  E_3 \otimes \Bone_3' & \otimes & \Bone_{12}'\\
  \Bone_3'\otimes E_3  & \otimes & \Bone_{12}'\\
  \Bone_9'             & \otimes & B\\
\end{array}
\right]
\]
and the configuration $A$ for the group-wise independence model is written as
\[
A = \left[
\begin{array}{c}
E_9 \otimes \Bone_{12}'\\
\Bone_9' \otimes E_{12}'
\end{array}
\right],
\]
where $E_n$ is the $n\times n$ identity matrix, $\Bone_n =
(1,\ldots,1)'$ is the $n\times 1$ column vector of $1$'s, 
$\otimes$ denotes the Kronecker product and
\[
 B = \left[\begin{array}{c}
111100000000\\
000011110000\\
100010001100\\
010001000011\\
001000101010\\
000100010101
\end{array}
\right].
\]
Note that the configuration $B$ is
the vertex-edge incidence matrix of
the $(2,2,2)$ complete multipartite graph.
Quadratic Gr\"obner bases of toric ideals arising from
complete multipartite graphs are studied in Ohsugi and Hibi (2000).

\begin{figure*}[htbp]
\begin{center}
\includegraphics[width=5.7cm]{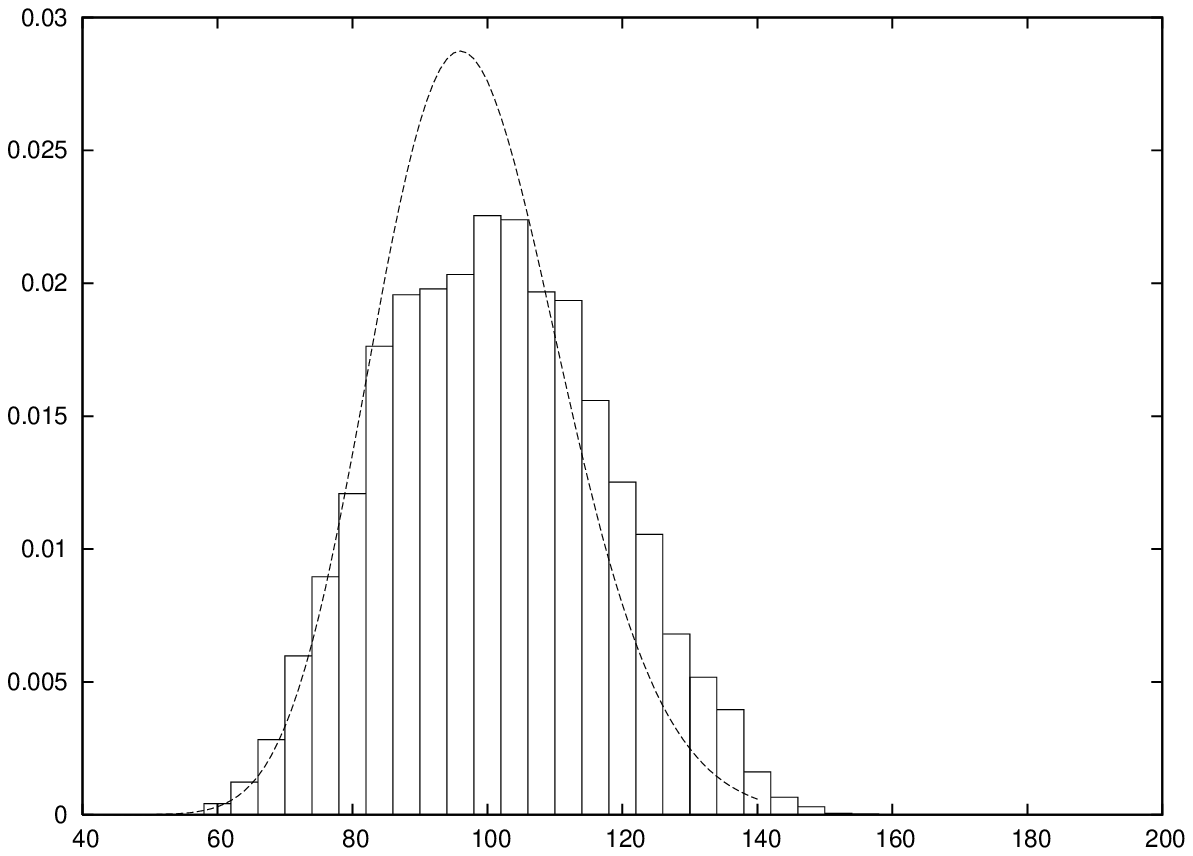}
\includegraphics[width=5.7cm]{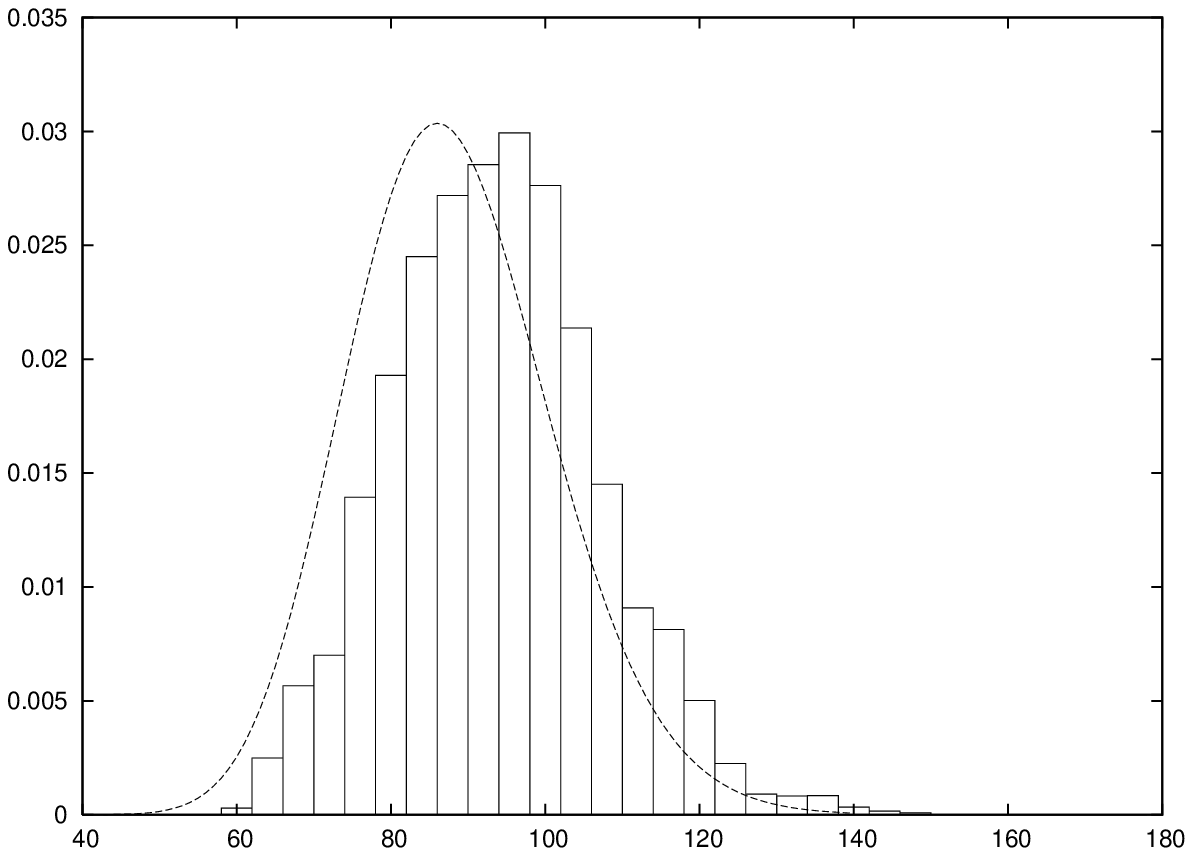}\\
\ \ Complete independence model ($df=98$) \quad Group-wise independence model ($df=88$)
\caption{Asymptotic and Monte Carlo sampling distributions of NCT
  data}
\label{fig:plot-NCT}
\end{center}
\end{figure*}

Given these configurations we can easily run a Markov chain as
discussed at the end of Section \ref{sec:GB}.
After $5,000,000$ burn-in steps, we construct $10,000$ Monte Carlo samples.
Figure \ref{fig:plot-NCT} show histograms of the Monte Carlo sampling 
generated from the exact conditional distribution of the Pearson goodness-of-fit 
$\chi^2$ statistics for the NCT data under the complete independence 
model and the group-wise independence model, respectively, along with the 
corresponding asymptotic distributions $\chi^2_{98}$ and $\chi^2_{88}$.

\subsection{The analysis of PTGDR (prostanoid DP receptor) diplotype frequencies data}
\label{subsec:diplo-data}
Next we give a numerical example of genome data.
Table \ref{tbl:PTGDR-data} shows diplotype frequencies on the three
loci, T-549C (locus 1), C-441T (locus 2) and T-197C (locus 3) in the
human genome 14q22.1, which is
given in Oguma {\it et al.} (2004).
Though the data is used for the genetic
association studies in Oguma {\it et al.} (2004), we simply consider fitting our
models. As an example, we only consider the diplotype data of patients
in the population of blacks ($n = 79$). 
\begin{table*}[htbp]
\begin{center}
\caption{PTGDR diplotype frequencies among patients and
 controls in each population. (The order of the SNPs in the haplotype is
 T-549C, C-441T and T-197C.)}
\label{tbl:PTGDR-data}
\begin{tabular}{ccccc}\hline
Diplotype & \multicolumn{2}{c}{Whites} & \multicolumn{2}{c}{Blacks}\\
& \multicolumn{1}{c}{Controls}& \multicolumn{1}{c}{Patients}
& \multicolumn{1}{c}{Controls}& \multicolumn{1}{c}{Patients}\\ \hline
CCT/CCT & 16 & 78 & 7 & 10\\
CCT/TTT & 27 & 106 & 12 & 27\\
CCT/TCT & 48 & 93 & 4 & 12\\
CCT/CCC & 17 & 45 & 3 & 9\\
TTT/TTT & 9 & 43 & 2 & 7\\
TTT/TCT & 34 & 60 & 8 & 6\\
TTT/CCC & 4 & 28 & 1 & 6\\
TCT/TCT & 11 & 20 & 7 & 0\\
TCT/CCC & 6 & 35 & 1 & 2\\
CCC/CCC & 1 & 8 & 0 & 0\\ \hline
\end{tabular}
\end{center}
\end{table*}

\begin{figure*}[htbp]
\begin{center}
\includegraphics[width=5.8cm]{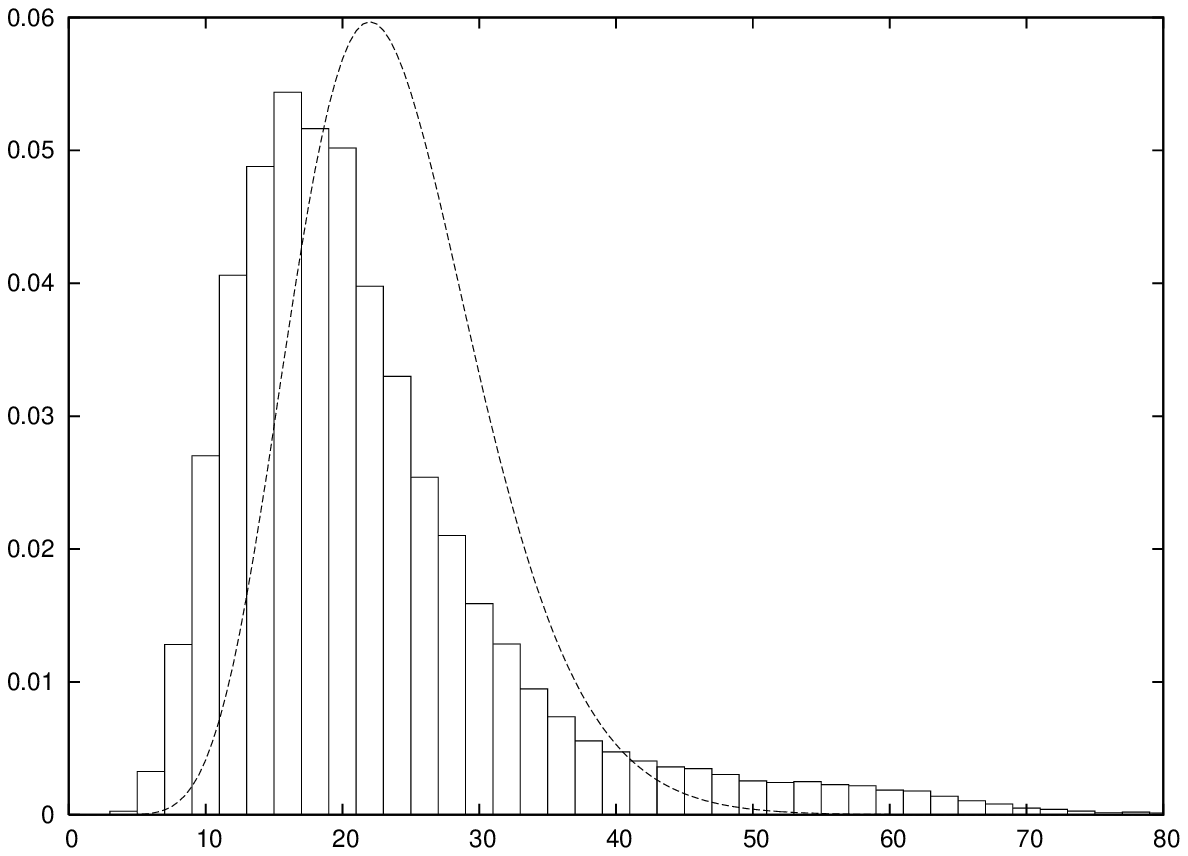}
\includegraphics[width=5.8cm]{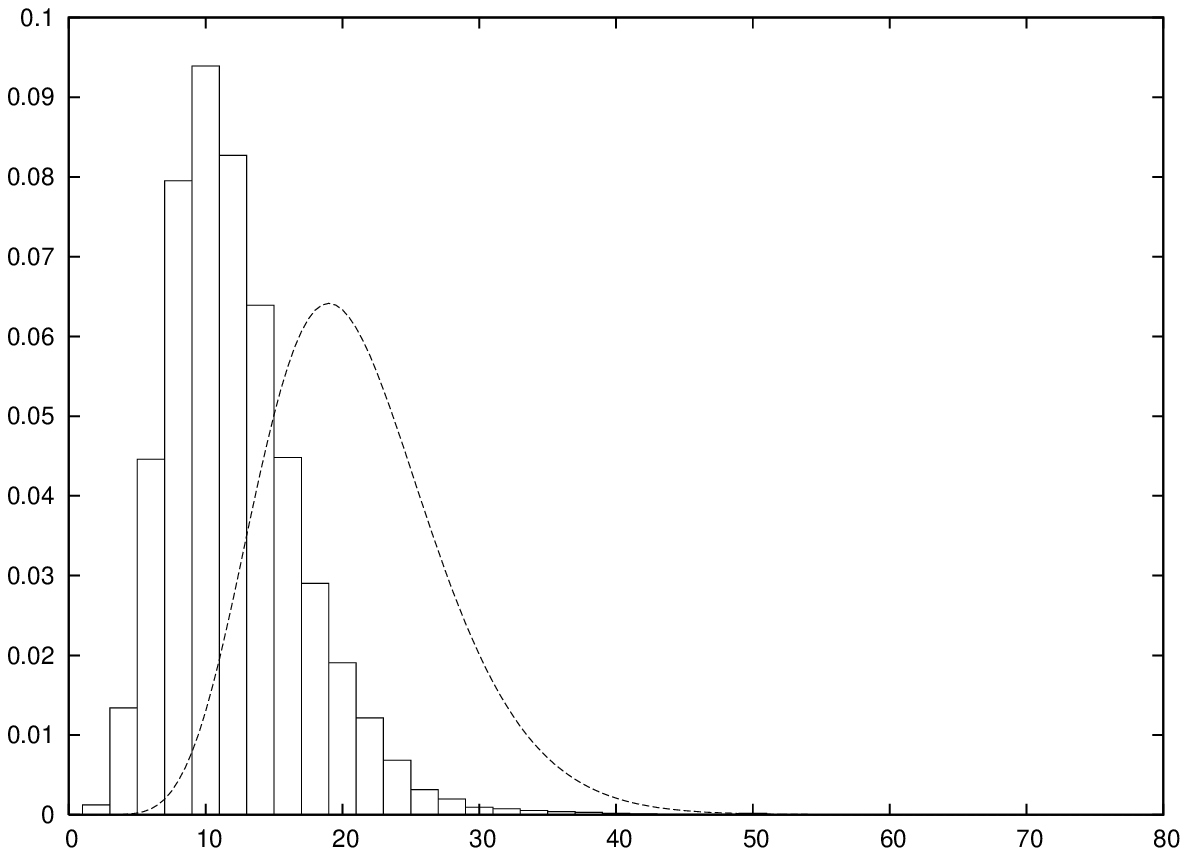}\\
Hardy-Weinberg model ($df=24$)\qquad\   Genotype-wise independence model ($df=21$)
\caption{Asymptotic and Monte Carlo sampling distributions of PTGDR genotype frequency data} 
\label{fig:plot-genotype}
\end{center}
\end{figure*}

First we consider the analysis of genotype frequency data. Though Table \ref{tbl:PTGDR-data}
is diplotype frequency data, here we ignore the information on the haplotypes and
simply treat it as a genotype frequency data. 
Since $\numgroup = 3$ and $m_1 = m_2 = m_3 = 2$, 
there are $3^3 = 27$ distinct set of genotypes, i.e., $|{\cal I}| = 27$, while only $8$ distinct
haplotypes appear in Table \ref{tbl:PTGDR-data}. 
Table \ref{tbl:genotype-PTGDR-data} is the set of genotype frequencies
of patients in the population of blacks.
\begin{table*}[htbp]
\begin{center}
\caption{The genotype frequencies for patients among blacks of PTGDR data}
\label{tbl:genotype-PTGDR-data}
\begin{tabular}{|c|c|c|c|c|c|c|c|c|c|c|}\hline
\multicolumn{2}{|c|}{locus 3} & \multicolumn{3}{|c|}{CC} & \multicolumn{3}{|c|}{CT} & \multicolumn{3}{|c|}{TT}\\ \hline
\multicolumn{2}{|c|}{locus 2} & CC & CT & TT             & CC & CT & TT             & CC & CT & TT\\ \hline
locus 1 & CC                  & 0  & 0  & 0              & 9  & 0  & 0              & 10 & 0  & 0 \\ \cline{2-11}
        & CT                  & 0  & 0  & 0              & 2  & 6  & 0              & 12 & 27 & 0 \\ \cline{2-11}
        & TT                  & 0  & 0  & 0              & 0  & 0  & 0              & 0  & 6  & 7 \\ \hline
\end{tabular}
\end{center}
\end{table*}
Under the genotype-wise independence model
(\ref{eq:genotype-indepedent}), the sufficient statistic is the genotype
frequency data for each locus. 
On the other hand, under the Hardy-Weinberg model
(\ref{eq:genetype-HW}), the sufficient statistic is the allele frequency
data for each locus, and the genotype frequencies for each locus are estimated
by the Hardy-Weinberg law.
Accordingly, the maximum likelihood estimates for the combination of
the genotype frequencies are calculated as Table \ref{tbl:genotype-PTGDR-data-MLE}.
\begin{table*}[htbp]
\begin{center}
\caption{MLE for PTGDR genotype frequencies
of patients among blacks under the Hardy-Weinberg model (upper) and genotype-wise
 independence model (lower)}
\label{tbl:genotype-PTGDR-data-MLE}
\begin{tabular}{|c|c|c|c|c|c|c|c|c|c|c|}\hline
\multicolumn{2}{|c|}{locus 3} & \multicolumn{3}{|c|}{CC} & \multicolumn{3}{|c|}{CT} & \multicolumn{3}{|c|}{TT}\\ \hline
\multicolumn{2}{|c|}{locus 2} & CC & CT & TT             & CC & CT & TT             & CC & CT & TT\\ \hline
locus 1 & CC & 0.1169 & 0.1180 & 0.0298 & 1.939 & 1.958 & 0.4941 & 8.042 & 8.118 & 2.049 \\ 
        &    & 0 & 0 & 0 & 1.708 & 2.018 & 0.3623 & 6.229 & 7.361 & 1.321 \\ \cline{2-11}
        & CT & 0.2008 & 0.2027 & 0.0512 & 3.331 & 3.362 & 0.8486 & 13.81 & 13.94 & 3.519 \\ 
        &    & 0 & 0 & 0 & 4.225 & 4.993 & 0.8962 & 15.41 & 18.21 & 3.268 \\ \cline{2-11}
        & TT & 0.0862 & 0.0870 & 0.0220 & 1.430 & 1.444 & 0.3644 & 5.931  & 5.988  & 1.511 \\ 
        &    & 0 & 0 & 0 & 1.169 & 1.381 & 0.2479 & 4.262 & 5.037 & 0.9040 \\ \hline
\end{tabular}
\end{center}
\end{table*}
The configuration $A$ for the Hardy-Weinberg model is written as
\[
A = \left[
\begin{array}{ccc}
222222222 & 111111111 & 000000000\\
000000000 & 111111111 & 222222222\\
222111000 & 222111000 & 222111000\\
000111222 & 000111222 & 000111222\\
210210210 & 210210210 & 210210210\\
012012012 & 012012012 & 012012012
\end{array}
\right]
\]
and the configuration $A$ for the genotype-wise independence model is written as
\[
A = \left[
\begin{array}{c}
E_3 \otimes \Bone_3' \otimes \Bone_3'\\
\Bone_3' \otimes E_3 \otimes \Bone_3'\\
\Bone_3' \otimes \Bone_3' \otimes E_3'
\end{array}
\right].
\]
Since these two configurations are of the Segre-Veronese type, again we can easily 
perform MCMC sampling as discussed in Section \ref{sec:GB}.
After $100,000$ burn-in steps, we construct $10,000$ Monte Carlo samples.
Figure \ref{fig:plot-genotype} shows histograms of the Monte Carlo sampling 
generated from the exact conditional distribution of the Pearson goodness-of-fit 
$\chi^2$ statistics for the PTGDR genotype frequency data under the Hardy-Weinberg
model and the genotype-wise independence model, respectively, along with the 
corresponding asymptotic distributions $\chi^2_{24}$ and $\chi^2_{21}$.

From the Monte Carlo samples, we can also estimate the $P$-values for
each null model.
The values of the Pearson goodness-of-fit $\chi^2$ for the PTGDR
genotype frequency data of
Table \ref{tbl:genotype-PTGDR-data} are $\chi^2 = 88.26$ under the
Hardy-Weinberg models,
whereas $\chi^2 = 103.37$ under the genotype-wise independence
model. These values are
highly significant ($p < 0.01$ for both models), which implies the
susceptibility of 
the particular haplotypes. 

\begin{figure*}[htbp]
\begin{center}
\includegraphics[width=8cm]{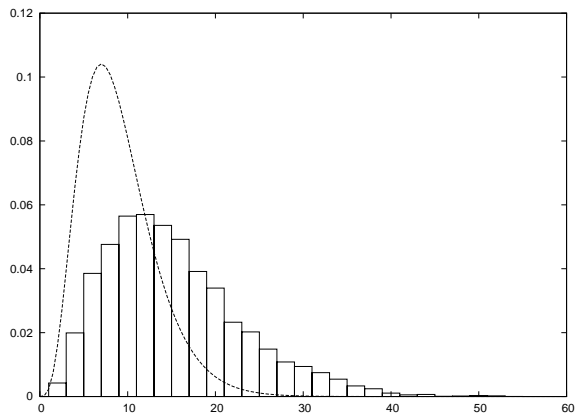}
\caption{Asymptotic and Monte Carlo sampling distributions of PTGDR diplotype frequency data
under the haplotype-wise Hardy-Weinberg model ($df = 9$).}
\label{fig:plot-diplotype}
\end{center}
\end{figure*}

Next we consider the analysis of the diplotype frequency data.
In this case of $\numgroup = 3$ and $m_1 = m_2 = m_3 = 2$, 
there are $2^3 = 8$
distinct haplotypes, and there are 
\[
 |{\cal I}| = 8 + {8\choose{2}} = 36
\]
distinct diplotypes, while there are only $4$ haplotypes
and $10$ diplotypes appear in Table \ref{tbl:PTGDR-data}. 
The numbers of each haplotype are calculated as the second column of Table
\ref{tbl:MLE-PTGDR-haplotype}.
Under the Hardy-Weinberg model, the haplotype frequencies are estimated proportionally
to the allele frequencies, which is shown as the third column of Table 
\ref{tbl:MLE-PTGDR-haplotype}.
\begin{table*}[htbp]
\begin{center}
\caption{Observed frequency and MLE under the 
Hardy-Weinberg model for PTGDR haplotype
 frequencies of patients among blacks.}
\label{tbl:MLE-PTGDR-haplotype}
\begin{tabular}{ccccccc}\cline{1-3}\cline{5-7}
Haplotype & observed & MLE under HW & & Haplotype & observed & MLE under HW \\ \cline{1-3}\cline{5-7}
CCC & 17 &   6.078  & & TCC & 0 &   5.220  \\
CCT & 68 &  50.410  & & TCT & 20 &  43.293  \\
CTC & 0 &   3.068   & & TTC & 0 &   2.635  \\
CTT & 0 &  25.445   & & TTT & 53 &  21.853  \\ \cline{1-3}\cline{5-7}
\end{tabular}
\end{center}
\end{table*}
The maximum likelihood estimates of the diplotype frequencies under the Hardy-Weinberg
model are calculated from the maximum likelihood estimates for each haplotype.
These values coincide with appropriate fractions of 
the values for the corresponding combination of the genotypes
in Table \ref{tbl:genotype-PTGDR-data-MLE}. For example, the MLE for the
diplotype CCT/CCT coincides with the MLE for the combination of the genotypes (CC,CC,TT) in
Table \ref{tbl:genotype-PTGDR-data-MLE}, whereas the MLE's for the diplotype CCC/TTT, CCT/TTC,
CTC/TCT, CTT/TCC coincide with the $\frac{1}{4}$ fraction of the 
MLE for the combination of the genotypes (CT,CT,CT), and so on.
Since we know that the Hardy-Weinberg model is highly
statistically rejected, it is natural to consider the haplotype-wise
Hardy-Weinberg model given in Section \ref{sec:haplotype}. 
Table \ref{tbl:MLE-PTGDR-diplotype}
shows the maximum likelihood estimates under the haplotype-wise Hardy-Weinberg model.
It should be noted that the MLE for the other diplotypes are all zeros.
\begin{table*}[htbp]
\begin{center}
\caption{MLE for PTGDR diplotype frequencies
of patients among blacks under the haplotype-wise Hardy-Weinberg model.}
\label{tbl:MLE-PTGDR-diplotype}
\begin{tabular}{ccccccc}\cline{1-3}\cline{5-7}
Diplotype & observed & MLE & & Diplotype & observed & MLE\\ \cline{1-3}\cline{5-7}
CCT/CCT & 10 &  14.6329 & & TTT/TCT &  6 &   6.7089\\
CCT/TTT & 27 &  22.8101 & & TTT/CCC &  6 &   5.7025\\
CCT/TCT & 12 &   8.6076 & & TCT/TCT &  0 &   1.2658\\
CCT/CCC &  9 &   7.3165 & & TCT/CCC &  2 &   2.1519\\
TTT/TTT &  7 &   8.8892 & & CCC/CCC &  0 &   0.9146\\ \cline{1-3}\cline{5-7}
\end{tabular}
\end{center}
\end{table*}
We perform the Markov chain Monte Carlo sampling for the haplotype-wise Hardy-Weinberg model.
The configuration $A$ for this model is written as
\[
A = \left[
\begin{array}{c}
200000001111111000000000000000000000\\
020000001000000111111000000000000000\\
002000000100000100000111110000000000\\
000200000010000010000100001111000000\\
000020000001000001000010001000111000\\
000002000000100000100001000100100110\\
000000200000010000010000100010010101\\
000000020000001000001000010001001011
\end{array}
\right],
\]
which is obviously of the Segre-Veronese type. We give 
a histogram of the Monte Carlo sampling 
generated from the exact conditional distribution of the Pearson goodness-of-fit 
$\chi^2$ statistics for the PTGDR diplotype frequency data under the haplotype-wise 
Hardy-Weinberg model, along with the 
corresponding asymptotic distributions $\chi^2_{9}$ in Figure \ref{fig:plot-diplotype}.

The $P$-value for this model is estimated as $0.8927$ with the estimated
standard deviation $0.0029$
(We also discard the first $100,000$ samples, and use a batching method
to obtain an estimate
of variance, see Hastings (1970) and Ripley (1987)). Note that the asymptotic
$P$-value based on $\chi_9^2$ is $0.6741$.

\section{Some discussions}
\label{sec:discussion}

In this paper we considered independence models in group-wise
selections, which can be described in terms of a Segre-Veronese
configuration.  We have shown that our framework can be applied to two
important examples in educational statistics and biostatistics.  We
expect that the methodology of the present paper finds applications in
many other fields.

In the NCT example, we assumed that the examinees choose the same
number $\totalitems$ of subjects.  We also assumed for simplicity that
the examinees choose either nothing or one 
subject from a  subgroup.  This restricts our analysis to some
subset of the examinees of NCT.  Actually the examinees make decisions
on how many subjects to take and modeling this decision making is
clearly of statistical interest.  Further complication arises from the
fact that the examinees can choose which scores to submit to
universities after taking NCT.  For example after obtaining scores of
three subjects on Science, an examinee can choose the best two scores for
submitting to a university.  In our subsequent paper 
(Aoki et al., 2007) we present a
generalization of  Segre-Veronese configurations
to cope with these complications.

It seems that the simplicity of the reduced Gr\"obner basis for the
Segre-Veronese configuration comes from the fact that the index set
$\cal J$ of the rows of $A$ can be ordered and the restriction on the
counts can be expressed in terms of one-dimensional intervals.  {}From
statistical viewpoint, ordering of the elements of the sufficient statistic in
group-wise selection seems to be somewhat artificial.  It is of
interest to look for other statistical models, where ordering of
the elements of the sufficient statistic is more natural and 
the Segre-Veronese configuration can be applied.

\begin{landscape}
\appendix

\renewcommand{\intextsep}{5pt}

\section{Tables of numbers of examinees in NCT in 2006}
\begin{table}[H]
\small
\begin{center}
\caption{Number of examinees who takes subjects on Social Studies}
\label{tbl:data-soc}
\begin{tabular}{|c|c|c|c|c|c|c|c|c|c|c|c|}\hline
 & \multicolumn{6}{|c|}{Geography and History}
 & \multicolumn{3}{|c|}{Civics} & \# total &
 \#  actual   \\ \cline{2-10}
 & WHA & WHB & JHA & JHB & GeoA & GeoB
& ContS & Ethics & P\&E & examinees &
 examinees  \\ \hline
1 subject & 496 & 29,108 & 1,456 & 54,577 & 1,347 & 27,152 
& 40,677 & 16,607 & 25,321 & 196,741 & 
196,741 \\ \hline
2 subjects & 1,028 & 61,132 & 3,386 & 90,427 & 5,039 & 83,828
 & 180,108 & 27,064 & 37,668 & 489,680 
& 244,840  
\\ \hline
Total & 1,524 & 90,240 & 4,842 & 145,004 & 6,386 & 110,980
& 220,785 & 43,671 & 62,989 & 686,421 & 441,581 \\ \hline
\end{tabular}
\end{center}
\end{table}
\begin{table}[H]
\begin{center}
\small
\caption{Number of examinees who selects two  subjects on Social
 Studies}
\label{tbl:data-soc2}
\begin{tabular}{|c|c|c|c|c|c|c|c|}\hline
\multicolumn{1}{|c|}{} & \multicolumn{6}{|c|}{Geography and History} & \\
 \cline{2-7}
\multicolumn{1}{|c|}{Civics} & WHA & WHB & JHA & JHB & GeoA
 & GeoB & Total\\  \hline
ContSoc  & 687 & 39,913 & 2,277 & 62,448 & 3,817 & 70,966 & 180,108\\
 \cline{1-8}
Ethics  & 130 & 10,966 & 409 & 10.482 & 405 &  4,672 & 27,064\\  
 \cline{1-8}
P\&E & 211 & 10253 & 700 & 17,497 & 817 & 8,190 & 37,668\\ \hline
\multicolumn{1}{|c|}{Total} & 1,028 & 61,132 & 3,386 & 90,427 & 5,039 &
 83,838 & 244,840\\ \hline
\end{tabular}
\end{center}
\end{table}
\begin{table}[H]
\small
\begin{center}
\caption{Number of examinees who takes subjects on Science}
\label{tbl:data-sci}
\begin{tabular}{|c|c|c|c|c|c|c|c|c|c|c|c|c|c|}\hline
 & \multicolumn{4}{|c|}{Science 1} & \multicolumn{3}{|c|}{Science 2}
& \multicolumn{4}{|c|}{Science 3} & \# total & \#actual 
\\ \cline{2-12}
 & CSciB   & BioI & ISci & BioIA & CSciA   & ChemI &
 ChemIA
& PhysI & EarthI & PhysIA & EarthIA & examinees & examinees 
\\\hline
1 subject & 2,558 & 80,385 & 511 & 1,314 & 1,569 & 19,616 & 717
& 14,397 & 10,788 & 289 & 236 & 132,380 &
132,380 
\\ \hline
2 subjects & 6,878 & 79,041 & 523 & 1,195 & 26,848 & 158,027 & 2,777
& 106,822 & 6,913 & 905 & 259 & 390,188 &
195,094 
\\\hline
3 subjects & 7,942 & 18,519 & 728 & 490 & 6,838 & 20,404 & 437
& 18,451 & 8,423 & 361 & 444 & 83,037 &  27,679 
\\ \hline
Total & 17,378 & 177,945 & 1,762 & 2,999 & 35,255 & 198,047 & 3,931
&139,670 & 26,124 & 1,555 & 939 & 605,605 & 355,153
\\\hline
\end{tabular}
\end{center}
\end{table}
\begin{table}[H]
\begin{center}
\caption{Number of examinees who selects two subjects on Science}
\label{tbl:data-sci2}
\small
\begin{tabular}{|c|c|c|c|c|c|c|c|c|}\hline
\multicolumn{2}{|c|}{} & \multicolumn{3}{|c|}{Science $2$} & 
\multicolumn{4}{|c|}{Science $3$}\\
 \cline{3-9}
\multicolumn{2}{|c|}{}            & CSciA   & ChemI & ChemIA & 
PhysI & EarthI     & PhysIA & EarthIA \\ \hline
Science 1 & CSciB  & 1,501 & 1,334 & 23 & 120 & 3,855 & 1 & 44\\
 \cline{2-9}
 & BioI  & 21,264 & 54,412 & 244 & 1,366 & 1,698 & 5 & 52\\ 
\cline{2-9}
& ISci & 147 & 165 & 50 & 43 & 92 & 5 & 21\\
\cline{2-9}
& BioIA & 128 & 212 & 715 & 16 & 33 & 29 & 62\\
\hline
Science 3& Physics & 3,243 & 101,100 & 934 & --- & --- & --- & ---\\
\cline{2-9}
 & EarthI     & 485 & 730 & 20 & --- & --- & --- & ---\\
\cline{2-9}
& PhysIA & 43 & 54 & 768 & --- & --- & --- & ---\\
\cline{2-9}
& EarthIA & 37 & 20 & 23 & --- & --- & --- & ---\\
\hline
\end{tabular}
\end{center}
\end{table}
\begin{table}[H]
\begin{center}
\caption{Number of examinees who selects three subjects on Science}
\label{tbl:data-sci3}
\small
\begin{tabular}{|c|c|c|c|c|c|c|c|c|c|c|c|c|c|c|}\hline
\multicolumn{2}{|c|}{Science $3$} & \multicolumn{3}{|c|}{PhysI} & 
\multicolumn{3}{|c|}{EarthI}
& \multicolumn{3}{|c|}{Physics IA} & 
\multicolumn{3}{|c|}{Earth science IA}\\
\hline
\multicolumn{2}{|c|}{Science $2$}  & CSciA   & ChemI & ChemIA & 
CSciA   & ChemI & ChemIA &
CSciA   & ChemI & ChemIA & 
CSciA   & ChemI & ChemIA\\ \hline 
Science $1$ & CSciB   & 1,155 & 5,152 & 17 & 1,201 & 317 & 7
&16 & 5 & 16 & 48 & 5 & 3\\ \cline{2-14}
& BioI  & 553 & 10,901 & 31 & 3,386 & 3,342 & 16 
&30 & 35 & 19 & 130 & 56 & 20 \\ \cline{2-14}
& ISci & 80 & 380 & 23 & 62 & 34 & 4 
& 32 & 13 & 27 & 48 & 14 & 11 \\ \cline{2-14}
& BioIA & 6 & 114 & 39 & 22 & 22 & 10 
&12 & 6 & 150 & 57 & 8 & 44 \\
\hline
\end{tabular}
\end{center}
\end{table}
\end{landscape}


\begin{thebibliography}{99}

\bibitem[Aoki et al.(2007)]{AHOT2}
Aoki, S., Hibi, T., Ohsugi, H.\ and Takemura, A. (2007).
Gr\"obner bases of nested configurations.
Submitted for publication.

\bibitem[Aoki and Takemura(2005)]{AT05}
Aoki, S.\  and  Takemura, A. (2005).
Markov chain Monte Carlo exact tests for incomplete two-way
contingency tables.
{\it Journal of Statistical Computation and Simulation}.
{\bf 75}, 787--812.

\bibitem[Aoki and Takemura(2006)]{AT06b}
Aoki, S.\ and Takemura, A. (2006).
Markov chain Monte Carlo tests for designed experiments.
{\tt arXiv:math/0611463v1}.  Submitted for publication.

\bibitem[Aoki and Takemura(2007)]{AT06}
Aoki, S.\ and Takemura, A. (2007).
Minimal invariant Markov basis for sampling contingency tables with
		       fixed marginals. {\it Annals of the Institute of
		       Statistical Mathematics}, To appear.




\bibitem{BES06}
Beerenwinkel, N., Eriksson, N.\ and Sturmfels, B. (2006).
Conjunctive Bayesian networks.
{\tt arXiv:math/0608417v3}. 
To appear in {\it Bernoulli}. 


\bibitem{BFH}
Bishop, Y.\ M.\ M., Fienberg, S.\ E.\  and  Holland
P.\ W. (1975). {\it Discrete Multivariate
Analysis: Theory and Practice}. The MIT Press, Cambridge, Massachusetts.

\bibitem[Crow (1988)]{C88}
Crow, J.\ E. (1988).
Eighty years ago: The beginnings of population genetics. {\it Genetics},
		       {\bf 119}, 473--476.

\bibitem{NeHi}
De Negri, E.\ and Hibi, T. (1997).
Gorenstein algebras of Veronese type.
{\it Journal of Algebra}, {\bf 193}, no. 2, 629--639. 


\bibitem[Diaconis and Sturmfels(1998)]{DS98}
Diaconis, P.\ and Sturmfels, B. (1998).  
Algebraic algorithms for sampling from conditional distributions. {\it
		       The Annals of Statistics}, {\bf 26}, 363--397.
\bibitem{Guo-Thompson}
Guo, S. and Thompson, E. (1992).
Performing the exact test of Hardy-Weinberg proportion for multiple
				 alleles. 
{\it Biometrics}, {\bf 48}, 361--372.


\bibitem[Hastings(1970)]{Has70}
Hastings, W.\ K. (1970). Monte Carlo sampling methods
  using Markov chains and their
  applications. {\it Biometrika}, {\bf
  57}, 97--109.


\bibitem{hibi87}
Hibi, T. (1987). Distributive lattices, affine semigroup rings and
algebras with straightening laws. {\it Advanced Studies in\ Pure
			   Mathematics}, {\bf
  11}, 93--109.

\bibitem{Huber-et-al-2006}
Huber, M., Chen, Y., Dinwoodie, I., Dobra, A.\ and Nicholas, M. (2006).
Monte Carlo algorithms for Hardy-Weinberg proportions.
{\it Biometrics}, {\bf 62}, 49--53.


\bibitem[Lehmann(1986)]{L86}
Lehmann, E.\ L.\ and Romano, J.\ P. (2005). {\it Testing Statistical Hypotheses}, 3rd  ed. Springer, New York.


\bibitem{NCUEE-H18}
National Center for University Entrance Examinations. (2006).  Booklet
on NCUEE.  Available from {\tt http://www.dnc.ac.jp/dnc\allowbreak{}/gaiyou\allowbreak{}/pdf\allowbreak{}/youran\_english\_H18\_HP.pdf}


\bibitem{PTGDR2004}
Oguma, T., Palmer, L. J., Birben, E., Sonna, L. A. Asano, K. and Lilly,
				 C. M. (2004).
Role of prostanoid DP receptor variants in susceptibility to asthma.
				 {\it The New England Journal of
				 Medicine}, {\bf 351}, 1752--1763.

\bibitem{OH3}
Ohsugi, H.\ and Hibi, T. (2000).
Compressed polytopes, initial ideals and complete multipartite graphs,
\textit{Illinois Journal of Mathematics}, \textbf{44}, 391--406. 


\bibitem[Ohsugi and Hibi(2005)]{OH05}
Ohsugi, H.\ and Hibi, T. (2005). 
Indispensable
 binomials of finite graphs. {\it Journal of Algebra and Its
  Applications}, {\bf 4}, 421--434.


\bibitem{survey}
Ohsugi, H.\ and Hibi, T. (2006).
Quadratic Gr\"obner bases 
arising from combinatorics.
{\it Integer Points in Polyhedra - Geometry, Number Theory, Representation
Theory, Algebra, Optimization, Statistics}, To appear.  

\bibitem{Ohsugi-Hibi-2007}
Ohsugi, H. and Hibi, T. (2007). 
Toric ideals arising from contingency tables. 
{\it in} ``Commutative Algebra and Combinatorics,"
Ramanujan Mathematical Society Lecture Notes Series, Number 4,
Ramanujan Mathematical Society,
India, in press.



\bibitem{PH2005}
Pachter, L.\ and Sturmfels, B.\ (2005).  {\it Algebraic Statistics for 
Computational Biology.}  Cambridge University Press, Cambridge.

\bibitem{Rapallo-2006}
Rapallo, F. (2006).
Markov bases and structural zeros.
{\it Journal of Symbolic Computation}, {\bf 41}, 164--172.

\bibitem[Ripley(1987)]{R87}
Ripley, B. D. (1987).
{\it Stochastic Simulation}. Wiley, New York.


\bibitem[Sturmfels(1995)]{St95}
Sturmfels, B. (1995). 
{\it Gr\"obner Bases and
    Convex   Polytopes}. American Mathematical Society, Providence, RI.

\bibitem[Takemura and Aoki(2004)]{TA04}
  Takemura, A.\ and Aoki, S. (2004).
Some characterizations of minimal
  Markov basis for sampling  from discrete conditional distributions.
  {\it Annals of the Institute of Statistical Mathematics}, {\bf 56}, 1--17.
\end{thebibliography}
\end{document}